\documentclass[11pt]{article}
\usepackage{indentfirst}
\usepackage{amsfonts}
\usepackage{url}
\usepackage{algorithm}
\usepackage{algpseudocode}

\newtheorem{teor}{Theorem}[section]
\newcommand{\intt}{\int\!\!\int_{\Omega}}
\author{E. Scheiber\thanks{\textit{Transilvania} University of Bra\c{s}ov, e-mail: scheiber@unitbv.ro}}
\title{On the Chebyshev approximation of a function with two variables}
\date{}
\begin{document}
\maketitle

\begin{abstract}
There is presented an approach to find an approximation 
polynomial of a function with two variables based on the two dimensional
discrete Fourier transform. The approximation
polynomial is expressed through  Chebyshev polynomials. 
An uniform convergence result is given.


\end{abstract}

\section{Introduction}
The purpose of the paper is to present some aspects about the construction of an approximation
polynomial for a function  with two variables. The approximation
polynomial is expressed through  Chebyshev polynomials.

Constructing an approximation
polynomial of a function with the corresponding applications is the subject of
the \textit{Chebfun} software,  presented in details in \cite{1}, \cite{200}. The \textit{Chebfun2} part
of the software deals with the construction of an approximation polynomial of a function with two
variables. According to \cite{2}, \cite{3}, to this end 
it is used a method based on Gaussian elimination as a low rank function approximation. 

\^{I}n \textit{Chebfun} the approximation polynomial
of a function with one variable is obtained using an one dimensional discrete Fourier transform. 
The approach of this paper will use a  two dimensional discrete Fourier transform. 

In spectral methods the Chebyshev polynomials are often used. 
The same form of the approximation polynomial are used in \cite{6}, \cite{101}, too.

After recalling of some formulas on the Fourier series for a function with two variables and on the two dimensional discrete
Fourier transform there is presented an algorithm to obtain an approximation polynomial of a function
with two variables and a convergence result. A Lagrange type interpolation problem 
for a function with two variables is studied.
Two applications are mentioned: a numerical integration formula on a rectangle and
a numerical computation of the partial derivatives.  

\section{Two dimensional Fourier series}

Let $f:\mathbb{R}^2\rightarrow R$ be a continuous periodical function in each variable with the period $2\pi.$ 
The  Fourier series attached to the function is \cite{100}, t.3
$$
f(x,y)\sim
\sum_{n,m=0}^{\infty}(a_{n,m}\cos{nx}\cos{my}+b_{n,m}\cos{nx}\sin{my}+
$$
$$
+c_{n,m}\sin{nx}\cos{my}+d_{n,m}\sin{nx}\sin{my})
$$
with the coefficients given by
$$
\begin{array}{lcl}
a_{0,0}=\frac{1}{4\pi^2}\intt f(x,y)\mathrm{d}x\mathrm{d}y &\quad & 
a_{n,m}=\frac{1}{\pi^2}\intt f(x,y)\cos{nx}\cos{my}\mathrm{d}x\mathrm{d}y\\
\\
a_{n,0}=\frac{1}{2\pi^2}\intt f(x,y)\cos{nx}\mathrm{d}x\mathrm{d}y &\quad & 
b_{n,m}=\frac{1}{\pi^2}\intt f(x,y)\cos{nx}\sin{my}\mathrm{d}x\mathrm{d}y\\
\\
a_{0,m}=\frac{1}{2\pi^2}\intt f(x,y)\cos{my}\mathrm{d}x\mathrm{d}y &\quad & 
c_{n,m}=\frac{1}{\pi^2}\intt f(x,y)\sin{nx}\cos{my}\mathrm{d}x\mathrm{d}y\\
\\
b_{0,m}=\frac{1}{2\pi^2}\intt f(x,y)\sin{my}\mathrm{d}x\mathrm{d}y &\quad & 
d_{n,m}=\frac{1}{\pi^2}\intt f(x,y)\sin{nx}\sin{my}\mathrm{d}x\mathrm{d}y\\
\\
c_{n,0}=\frac{1}{2\pi^2}\intt f(x,y)\sin{nx}\mathrm{d}x\mathrm{d}y &\quad & 
\end{array}
$$
where $\Omega=[0,2\pi]^2.$

The complex form of the Fourier series is 
$$
\sum_{n.m\in\mathbb{Z}}\gamma_{n,m}e^{inx+imy}
$$
with
$$
\begin{array}{lcl}
\gamma_{0,0}=a_{0,0} & \quad &\\
\\
\gamma_{n,0}=\frac{1}{2}(a_{n,0}-ic_{n,0}) & & \gamma_{-n,0}=\frac{1}{2}(a_{n,0}+ic_{n,0}) \\
\\
\gamma_{0,m}=\frac{1}{2}(a_{0,m}-ib_{0,m}) & & \gamma_{0,-m}=\frac{1}{2}(a_{0,m}+ib_{0,m}) \\
\\
\gamma_{n,m}=\frac{1}{4}(a_{n,m}-ib_{n,m}-ic_{n,m}-d_{n,m}) & & \gamma_{-n,m}=\frac{1}{4}(a_{n,m}-ib_{n,m}+ic_{n,m}+d_{n,m}) \\
\\
\gamma_{n,-m}=\frac{1}{4}(a_{n,m}+ib_{n,m}-ic_{n,m}+d_{n,m}) & & \gamma_{-n,-m}=\frac{1}{4}(a_{n,m}+ib_{n,m}+ic_{n,m}-d_{n,m}) \\
\end{array}
$$
or
\begin{equation}\label{cheb2101}
\gamma_{m,n}=\frac{1}{4\pi^2}\intt f(x,y)e^{-inx-imy}\mathrm{d}x\mathrm{d}y,\quad \forall n,m\in\mathbb{Z}.
\end{equation}

If the function is even in any variable  then the $b_{n,m}, c_{n,m}, d_{n,m}$ coefficients are all zero.

We shall suppose that the convergence conditions of the Fourier series to $f(x,y)$ are fulfilled (the function has bounded
first order partial derivatives in $\Omega$ and in a neighborhood of $(x,y)$ there exists $\frac{\partial^2 f}{\partial x\partial y},$ 
or $\frac{\partial^2 f}{\partial y\partial x},$ which is continuous in $(x,y),$ cf. \cite{100}, t.3, 697).

\section{Two dimensional discrete Fourier transform}

Let be the infinite matrix $(x_{k,j})_{k,j\in\mathbb{Z}}$ with the periodicity properties
$x_{k+p,j}=x_{k,j},\ x_{k,j+q}=x_{k,j},\ \forall k,j\in\mathbb{Z}.$ The discrete Fourier transform
construct another infinite matrix $(y_{r,s})_{r,s\in\mathbb{Z}}$ with an analog periodicity properties
defined by
$$
y_{r,s}=\sum_{k=0}^{p-1}\sum_{j=0}^{q-1}x_{k,j}e^{-i\frac{2\pi k r}{p}}e^{-i\frac{2\pi j s}{q}},
$$
for $r\in\{0,1,\ldots,p-1\}$ and $s\in\{0,1,\ldots,q-1\}.$

The complexity to compute the $p q$ numbers with the discrete fast Fourier transform algorithm is $pq \log_2{pq}.$ 

As an application, if the Fourier series coefficients (\ref{cheb2101}) are computed using the trapezoidal rule
for each of the iterated integrals then:
$$
\gamma_{n,m}=\frac{1}{4\pi^2}\intt f(x,y)e^{-inx-imy}\mathrm{d}x\mathrm{d}y=
\frac{1}{4\pi^2}\int_0^{2\pi}e^{-inx}\left(\int_0^{2\pi}f(x,y)e^{-imy}\mathrm{d}y\right)\mathrm{d}x\approx
$$
$$
\approx \frac{1}{4\pi^2}\int_0^{2\pi}e^{-inx}\left(
\frac{2\pi}{q}\sum_{j=0}^{q-1}f(x,\frac{2\pi j}{q})e^{-i m\frac{2\pi j}{q}}
\right)\mathrm{d}x=
$$
$$
=\frac{1}{2\pi q}\sum_{j=0}^{q-1}e^{-i m\frac{2\pi j}{q}}\int_0^{2\pi}f(x,\frac{2\pi j}{q})e^{-inx}\mathrm{d}x\approx
\frac{1}{p q}\sum_{k=0}^{p-1}\sum_{j=0}^{q-1}f(\frac{2\pi k}{p},\frac{2\pi j}{q})e^{-i n\frac{2\pi k}{p}}e^{-i m\frac{2\pi j}{q}}.
$$
Thus the Fourier coefficients $(\gamma_{n,m})$ may be computed applying the discrete Fourier transform to
$\left(f(\frac{2\pi k}{p},\frac{2\pi j}{q})\right)_{k\in\{0,1,\ldots,p-1\}, j\in\{0,1,\ldots,q-1\}}.$

\section{The Chebyshev series}

Considering a continuous two real variables function $f(x,y),\ x,y\in[-1,1]$ there is attached the Chebyshev series
\begin{equation}\label{cheb2103}
f(x,y)\sim \sum_{n,m=0}^{\infty}\alpha_{n,m}T_n(x)T_m(y)
\end{equation}
where
$$
\begin{array}{lcl}
\alpha_{0,0}=\frac{1}{\pi^2}\int_{-1}^1\!\!\int_{-1}^1\frac{f(x,y)}{\sqrt{1-x^2}\sqrt{1-y^2}}\mathrm{d}x\mathrm{d}y & \quad &
\alpha_{n,0}=\frac{2}{\pi^2}\int_{-1}^1\!\!\int_{-1}^1\frac{f(x,y)T_n(x)}{\sqrt{1-x^2}\sqrt{1-y^2}}\mathrm{d}x\mathrm{d}y\\
\\
\alpha_{0,m}=\frac{2}{\pi^2}\int_{-1}^1\!\!\int_{-1}^1\frac{f(x,y)T_m(y)}{\sqrt{1-x^2}\sqrt{1-y^2}}\mathrm{d}x\mathrm{d}y & &
\alpha_{n,m}=\frac{4}{\pi^2}\int_{-1}^1\!\!\int_{-1}^1\frac{f(x,y)T_n(x)T_m(y)}{\sqrt{1-x^2}\sqrt{1-y^2}}\mathrm{d}x\mathrm{d}y 
\end{array}
$$
Changing $x=\cos{t},y=\cos{s},$ the coefficient $\alpha_{n,m}$ will be
\begin{equation}\label{cheb2104}
\alpha_{n,m}=\frac{4}{\pi^2}\int_0^{\pi}\!\!\int_0^{\pi}f(\cos{t},\cos{s})\cos{nt}\cos{ms}\mathrm{d}t\mathrm{d}s=
\end{equation}
$$
=\frac{1}{\pi^2}\intt f(\cos{t},\cos{s})\cos{nt}\cos{ms}\mathrm{d}t\mathrm{d}s.
$$
Analogous formulas may be obtained for $\alpha_{0,0},\alpha_{n,0}$ and $\alpha_{0,m},$ too.
Thus the coefficients of the Chebyshev series are the coefficients of the Fourier coefficients of the function
$\varphi(t,s)= f(\cos{t},\cos{s}).$

If the function $f$ has second order derivatives then the Fourier series attached to $\varphi$ converges to $\varphi$ 
and consequently
\begin{equation}\label{cheb2102}
f(x,y)=\sum_{n,m=0}^{\infty}\alpha_{n,m}T_n(x)T_m(y),\quad x,y\in [-1,1].
\end{equation}

The polynomial
$$
f_{n,m}(x,y)=\sum_{k=0}^n\sum_{j=0}^m\alpha_{k,j}T_k(x)T_j(y)
$$
 is called 
the Chebyshev approximation polynomial of the function   $f(x,y)$ in the square $[-1,1]^2.$

The parameters $n,m$ are determined adaptively to satisfy the inequalities $|\alpha_{k,j}|<tol (=10^{-15},$
 machine precision), for $k>n$ and $j>m.$
This is the goal of the algorithm \ref{algcoefflag}. The coefficients whose absolute value are less then $tol$ are
eliminated and the remained coefficients are stored as a sparse matrix.

The $f_{n,m}(x,y)$ polynomial may be obtained with the least square method
as the solution of the optimization problem
$$
\min_{\lambda_{k,j}}\int_{-1}^1\!\!\int_{-1}^1\frac{1}{\sqrt{1-x^2}\sqrt{1-y^2}}\left( f(x,y)-
\sum_{k=0}^n\sum_{j=0}^m\lambda_{k,j}T_k(x)T_j(y)\right)^2\mathrm{d}x\mathrm{d}y.
$$

Due to the Parseval equality 
$$
\alpha_{0,0}^2+\frac{1}{2}\sum_{n=0}^{\infty}\alpha_{n,0}^2+\frac{1}{2}\sum_{m=0}^{\infty}\alpha_{0,m}^2+
\frac{1}{4}\sum_{n,m=1}^{\infty}\alpha_{n,m}^2=\frac{1}{\pi^2}\int_{-1}^1\int_{-1}^1\frac{f^2(x,y)}{\sqrt{1-x^2}\sqrt{1-y^2}}
\mathrm{d}x\mathrm{d}y
$$
the quality of the approximation polynomial may be evaluated by 
\begin{equation}\label{che2106}
\frac{1}{\pi^2}\int_{-1}^1\int_{-1}^1\frac{f^2(x,y)}{\sqrt{1-x^2}\sqrt{1-y^2}}-
\left(\alpha_{0,0}^2+\frac{1}{2}\sum_{k=0}^{n}\alpha_{k,0}^2+\frac{1}{2}\sum_{j=0}^{m}\alpha_{0,j}^2+
\frac{1}{4}\sum_{k=1}^n\sum_{j=1}^{m}\alpha_{k,j}^2\right).
\end{equation}

\begin{algorithm}
\caption{Algorithm to compute the Chebyshev approximation polynomial} \label{algcoefflag}
\begin{algorithmic}[1]
\Procedure {chebfun2}{f}
\State $n\leftarrow 8$
\State $tol\leftarrow 10^{-15}$
\State $sw\leftarrow true$
\While {$sw$} \Comment{The approximation polynomial is determined adaptively}
\State $m\leftarrow 2n$
\State $x,y \leftarrow \cos{\frac{2k\pi}{m}},\ k=0:m-1$
\State $z \leftarrow f(x,y)$
\State $g\leftarrow FFT(z)/m^2$
\State $a\leftarrow  4\Re g(1:n,1:n)$
\State $a(1,1) \leftarrow a(1,1))/4$
\State $a(1,2:n)\leftarrow a(1,2:n)/2$
\State $a(2:n,1)\leftarrow a(2:n,1)/2$
\If {$|a(i-1:i,1:n)|<tol\  \&\  |a(1:n,i-1:i)|<tol$} 
\State $sw\leftarrow false$
\Else
\State $n\leftarrow 2n$
\EndIf
\EndWhile
\For{$i=1:n$}\Comment{Removal of negligible coefficients}
\For{$j=1:n$}
\If {$|a(i,j)|<tol$}
\State{$a(i,j) \leftarrow 0$}
\EndIf
\EndFor
\EndFor\\
\hspace*{0.5cm}
\Return {$a$}
\EndProcedure
\end{algorithmic}
\end{algorithm}

The value of the polynomial $f_{n,m}$ in a point $(x,y)$ may be computed adapting the Clenshaw algorithm, \cite{101},
but we find that the evaluation of the expression $f_{n,m}(x,y)=V'_n(x)A_{n,m}V_m(y),$ where
$V_{\nu}(s)=(T_0(s),T_1(s),\ldots,T_{\nu}(s))'$ and $A=(a_{k,j})_{k=0:n,j=0:m}$ is more
efficient within a matrix oriented software. $V'$ denotes the transpose of the vector $V.$ The complexity order of both
algorithms is $O(nm).$

\section{The Chebyshev series of partial derivatives}

We assume that the function $f(x,y)$ has first order continuous partial derivatives and the series
(\ref{cheb2102}), there is required to find the coefficients $(b_{n,m})_{n,m\in\mathbb{N}}$ 
such that 
\begin{equation}\label{cheb2110}
\frac{\partial f(x,y)}{\partial x}=\sum_{n,m=0}^{\infty}b_{n,m}T_n(x)T_m(y).
\end{equation}
Using the equalities
$$
T_1'(x)=T_0(x),\qquad T_2'(x)=T_1(x)
$$
and
$$
\frac{1}{2}\left(\frac{T'_{n+1}(x)}{n+1}-\frac{T'_{n-1}(x)}{n-1}\right)=T_n(x),\qquad n>1.
$$
(\ref{cheb2110}) may be written as
$$
\frac{\partial f(x,y)}{\partial x}=\sum_{m=0}^{\infty}\left(\sum_{n=0}^{\infty}b_{n,m}T_n(x)\right)T_m(y)=
$$
$$
=\sum_{m=0}^{\infty}\left(b_{0,m}T'_1(x)+\frac{b_{1,m}}{2}\frac{T'_2(x)}{2}+\sum_{k=2}^{\infty}\frac{b_{k,m}}{2}
(\frac{T'_{k+1}(x)}{k+1}-\frac{T'_{k-1}(x)}{k-1})\right)T_m(y)=
$$
$$
=\sum_{m=0}^{\infty}\left((b_{0,m}-\frac{b_{2,m}}{2})T'_1(x)+\sum_{k=2}^{\infty}\frac{1}{2k}(b_{k-1,m}-b_{k+1,m})T'_k(x)\right)T_m(y)=
$$
$$
=\sum_{m=0}^{\infty}\left(\sum_{k=1}^{\infty}\alpha_{k,m}T'_k(x)\right)T_m(y).
$$
Identifying the coefficients of $T'_k(x)$ there is obtained the linear algebraic system
\begin{equation}\label{cheb2111}
\left\{ \begin{array}{lcl}
b_{0,m}-\frac{b_{2,m}}{2} & = & \alpha_{1,m} \\
\frac{1}{2k}(b_{k-1,m}-b_{k+1,m}) & = & \alpha_{k,m},\qquad k\ge 2, \qquad m\in\mathbb{N}.
\end{array}\right.
\end{equation}
Summing the above equalities for $k=n+1,n+3,n+5,\ldots$ it results
$$
b_{n,m}=(n+1)a_{n+1,m}+(n+3)a_{n+3,m}+(n+5)a_{n+5,m}+\ldots\qquad \forall\ n,m\in\mathbb{N}.
$$

In the same way it is deduced that the coefficients of the series 
$$
\frac{\partial f(x,y)}{\partial y}=\sum_{n,m=0}^{\infty}c_{n,m}T_n(x)T_m(y)
$$
satisfies the relations
\begin{equation}\label{cheb2112}
\left\{ \begin{array}{lcl}
c_{n,0}-\frac{c_{n,2}}{2} & = & \alpha_{n,1} \\
\frac{1}{2j}(c_{n,j-1}-c_{n,j+1}) & = & \alpha_{n,j},\qquad j\ge 2, \qquad m\in\mathbb{N}.
\end{array}\right.
\end{equation}
Let be 
\begin{equation}\label{cheb2114}
\frac{\partial f^2(x,y)}{\partial x\partial y}=\frac{\partial f^2(x,y)}{\partial y\partial x}=
\sum_{n,m=0}^{\infty}d_{n,m}T_n(x)T_m(y).
\end{equation}
Because $\frac{\partial f^2(x,y)}{\partial x\partial y}=\frac{\partial}{\partial y}(\frac{\partial f(x,y)}{\partial x})$
applying (\ref{cheb2112}) it results
$$
\frac{1}{2j}(d_{k,j-1}-d_{k,j+1})=b_{k,j}
$$
and consequently, for $k,j>1,$
\begin{equation}\label{cheb2113}
\frac{1}{4kj}(d_{k-1,j-1}-d_{k-1,j+1}-d_{k+1,j-1}+d_{k+1,j+1})=\frac{1}{2k}(b_{k-1,j}-b_{k+1,j})=\alpha_{k,j}.
\end{equation}

Denoting $\triangle_{k,j}=d_{k-1,j-1}-d_{k-1,j+1}-d_{k+1,j-1}+d_{k+1,j+1}),$ from the Parseval equality
corresponding to (\ref{cheb2114}) it results that
$$
\sum_{k=2}^{\infty}\sum_{j=2}^{\infty}\triangle_{k,j}^2\le16M_{1,1}^2,
$$
where $M_{1,1}\ge \max\{|\frac{\partial^2 f}{\partial x\partial y}(x,y)|: \ x,y\in[-1,1]\}.$

\section{The convergence of the Chebyshev series}

Using the techniques presented in \cite{4} and \cite{1}, in some hypotheses 
it may be proven that the convergence in (\ref{cheb2102}) is uniform in $[-1,1]^2,$
for $n,m\rightarrow\infty.$ First we state 

\begin{teor}\label{t1}
If the function $f$ has second order  continuous derivatives then 
\begin{equation}\label{cheb2105}
|\alpha_{n,0}|\le \frac{2 M_{2,0}}{(n-1)^2} \quad \mbox{and}\quad |\alpha_{n,1}|\le  \frac{8 M_{2,0}}{\pi(n-1)^2}, \quad n>1,
\end{equation}
\begin{equation}\label{cheb2106}
|\alpha_{0,m}|\le \frac{2 M_{0,2}}{(m-1)^2} \quad \mbox{and}\quad |\alpha_{1,m}|\le  \frac{8 M_{0,2}}{\pi(m-1)^2}, \quad m>1,
\end{equation}
where $M_{2,0}\ge \max\{|\frac{\partial^2 f}{\partial x^2}(x,y)|: \ x,y\in[-1,1]\},$
$M_{0,2}\ge \max\{|\frac{\partial^2 f}{\partial y^2}(x,y)|: \ x,y\in[-1,1]\}.$
\end{teor}

\vspace{0.3cm}\noindent
\textbf{Proof.} The coefficient $\alpha_{n,0}$ may be written as
$$
\alpha_{n,0}=\frac{2}{\pi^2}\int_0^{\pi}\left(\int_0^{\pi}f(\cos{t},\cos{s})\cos{ns}\mathrm{d}t\right)\mathrm{d}s.
$$
Two partial integrations are performed in the internal integral
$$
\int_0^{\pi}f(\cos{t},\cos{s})\cos{nt}\mathrm{d}t=
$$
$$
=\frac{1}{2n}\int_0^{\pi} \frac{\partial^2 f}{\partial y^2}(\cos{t},\cos{s})\sin{t}
\left(\frac{\sin{(n-1)t}}{n-1}-\frac{\sin{(n+1)t}}{n+1}\right)\mathrm{d}t.
$$
It results that
\begin{equation}\label{cheb2115}
|\int_0^{\pi}f(\cos{t},\cos{s})\cos{nt}\mathrm{d}t|\le \frac{\pi}{2n}M_{2,0}\left(\frac{1}{n-1}+\frac{1}{n+1}\right)\le
\frac{\pi M_{2,0}}{(n-1)^2}
\end{equation}
and consequently $|\alpha_{n,0}|\le \frac{2 M_{2,0}}{(n-1)^2}.$ 

Using (\ref{cheb2115}) in  
$\alpha_{n,1}=\frac{4}{\pi^2}\int_0^{\pi}\left(\int_0^{\pi}f(\cos{t},\cos{s})\cos{nt}\mathrm{d}t\right)\cos{s}\mathrm{d}s$
it results
$$
|\alpha_{n,1}|\le \frac{4}{\pi}\frac{M_{2,0}}{(n-1)^2}\int_0^{\pi}|\cos{s}|\mathrm{d}s=\frac{8M_{2,0}}{\pi(n-1)^2}.
$$

The proof of (\ref{cheb2106}) is similar.
\qquad \rule{5pt}{5pt}

\begin{teor}\label{cheb2lag14}
 If the function $f$ has second order continuous partial derivatives then $\lim_{n,m\rightarrow\infty}f_{n,m}= f$
uniformly in $[-1,1]^2.$
\end{teor}

\vspace{0.3cm}\noindent
\textbf{Proof.} From 
$$
f(x,y)=\sum_{n,m=0}^{\infty}\alpha_{n,m}T_n(x)T_m(y)\quad \mbox{and}\quad
f_{n,m}(x,y)=\sum_{k=0}^n\sum_{j=0}^m\alpha_{k,j}T_k(x)T_j(y)
$$
 it results
$$
f(x,y)-f_{n,m}(x,y)=\sum_{k=0}^n\sum_{j=m+1}^{\infty}\alpha_{k,j}T_k(x)T_j(y)+
\sum_{k=n+1}^{\infty}\sum_{j=0}^{\infty}\alpha_{k,j}T_k(x)T_j(y).
$$
Then
\begin{equation}\label{cheb2lag13}
|f(x,y)-f_{n,m}(x,y)|
\le \sum_{k=0}^n\sum_{j=m+1}^{\infty}|\alpha_{k,j}|+ 
\sum_{k=n+1}^{\infty}\sum_{j=0}^{\infty}|\alpha_{k,j}|=
\end{equation}
$$
=\sum_{j=m+1}^{\infty}|\alpha_{0,j}|+
\sum_{j=m+1}^{\infty}|\alpha_{1,j}|+
 \sum_{k=2}^n\sum_{j=m+1}^{\infty}|\alpha_{k,j}|+
 $$
 $$
+ \sum_{k=n+1}^{\infty}|\alpha_{k,0}|+ \sum_{k=n+1}^{\infty}|\alpha_{k,1}|+
\sum_{k=n+1}^{\infty}\sum_{j=2}^{\infty}|\alpha_{k,j}|
$$
and using the Cauchy-Buniakowsky-Schwarz inequality it follows
$$
(f(x,y)-f_{n,m}(x,y))^2\le
$$
$$
\le 6\left(
\left(\sum_{j=m+1}^{\infty}|\alpha_{0,j}|\right)^2+
\left(\sum_{j=m+1}^{\infty}|\alpha_{1,j}|\right)^2+
 \left(\sum_{k=2}^n\sum_{j=m+1}^{\infty}|\alpha_{k,j}|\right)^2+\right.
 $$
 $$
\left.+\left( \sum_{k=n+1}^{\infty}|\alpha_{k,0}|\right)^2+
\left( \sum_{k=n+1}^{\infty}|\alpha_{k,1}|\right)^2+
\left(\sum_{k=n+1}^{\infty}\sum_{j=2}^{\infty}|\alpha_{k,j}|\right)^2
\right)
$$
The following inequality holds $\sum_{i=\nu+1}^{\infty}\frac{1}{i^2}<\int_{\nu}^{\infty}\frac{\mathrm{d}x}{x}=\frac{1}{\nu}.$

Using the results of Theorem \ref{t1}, the first, second, fourth and the fifth expression
are increased by 
\begin{eqnarray*}
\sum_{j=m+1}^{\infty}|\alpha_{0,j}|& \le & 2M_{0,2}\sum_{j=m}^{\infty}\frac{1}{j^2}<\frac{2M_{0,2}}{m-1}  \\
\sum_{j=m+1}^{\infty}|\alpha_{1,j}|& \le & \frac{8M_{0,2}}{\pi}\sum_{j=m}^{\infty}\frac{1}{j^2}<\frac{8M_{0,2}}{\pi(m-1)}<
\frac{4M_{0,2}}{m-1}\\
 \sum_{k=n+1}^{\infty}|\alpha_{k,0}|& \le & 2M_{2,0}\sum_{k=n}^{\infty}\frac{1}{k^2}<\frac{2M_{2,0}}{n-1} \\
 \sum_{k=n+1}^{\infty}|\alpha_{k,1}|& \le & \frac{8M_{2,0}}{\pi}\sum_{k=n}^{\infty}\frac{1}{k^2}<\frac{8M_{2,0}}{\pi(n-1)}<
 \frac{4M_{2,0}}{n-1}.
\end{eqnarray*}
For the third and sixth expression we shell use (\ref{cheb2113}) and then  the Cauchy-Buniakowsky-Schwarz' inequality 
$$
\left(\sum_{k=2}^n\sum_{j=m+1}^{\infty}|\alpha_{k,j}|\right)^2=
\left(\sum_{k=2}^n\sum_{j=m+1}^{\infty}\frac{|\triangle_{k,j}|}{4kj}\right)^2\le
\frac{1}{16}\left(\sum_{k=2}^n\sum_{j=m+1}^{\infty}\triangle_{k,j}^2\right)
\left( \sum_{k=2}^n\sum_{j=m+1}^{\infty}\frac{1}{k^2j^2}\right)\le
$$
$$
\le M_{1,1}^2 \sum_{k=2}^n\frac{1}{k^2}\sum_{j=m+1}^{\infty}\frac{1}{j^2}\le
 \frac{\pi^2M_{1,1}^2}{6m},
$$
and respectively
$$
 \left(\sum_{k=n+1}^{\infty}\sum_{j=2}^{\infty}|\alpha_{k,j}|\right)^2\le
 \left(\sum_{k=n+1}^{\infty}\sum_{j=2}^{\infty}\frac{|\triangle_{k,j}|}{4kj}\right)^2\le
 M_{1,1}^2\left(\sum_{k=n+1}^{\infty}\frac{1}{k^2}\right)\left(\sum_{j=2}^{\infty}\frac{1}{j^2}\right)\le
 \frac{\pi^2M_{1,1}^2}{6n}
$$

Consequently
$$
|f(x,y)-f_{n,m}(x,y)|\le \sqrt{6}\left(\frac{20M_{0,2}^2}{(m-1)^2}+\frac{20M_{2,0}^2}{(n-1)^2}+ 
 \frac{\pi^2M_{1,1}^2}{6m}+ \frac{\pi^2M_{1,1}^2}{6n} \right)^{\frac{1}{2}} \rightarrow 0,
$$
when $m,n\rightarrow\infty.$ \rule{5pt}{5pt}

\section{The Lagrange interpolation polynomial}

For any grids $-1\le x_0<x_1<\ldots<x_n\le 1, \ -1\le y_0<y_1<\ldots<y_m\le 1$  and any 
$f:[-1,1]^2\rightarrow \mathbb{R}$ the expression of the Lagrange interpolation polynomial
is
$$
L_{n,m}(x,y)=\sum_{i=0}^n\sum_{j=0}^mf(x_i,y_j)lx_i(x)ly_j(y)
$$
where
$$
lx_i(x)=\prod_{k=0,k\not=i}^n\frac{x-x_k}{x_i-x_k}, \quad \mbox{and}\quad ly_j(x)=\prod_{l=0,l\not=j}^m\frac{y-y_l}{y_j-y_l}.
$$
This polynomial  satisfies the interpolation restrictions
$$
L_{n,m}(x_k,y_l)=f(x_k,y_l),\quad\forall k\in\{0,1,\ldots,n\},\ \mbox{and}\ \forall l\in\{0,1,\ldots,m\}.
$$
In the set of $(n,m)$ degree polynomials there exists a unique interpolation polynomial.

If $x_i=\cos{\frac{i\pi}{n}},\ i\in\{0,1,\ldots,n\}$ and $y_j=\cos{\frac{j\pi}{m}},\ j\in\{0,1,\ldots,m\}$
then 
using the \textit{discrete orthogonality} relations, \cite{101},
$$
\sum_{k=0}^n\gamma_kT_p(x_k)T_q(x_k)=\left\{\begin{array}{lcl}
0 &\mbox{if} & p\not=q\\
\frac{n}{2} &\mbox{if}& p=q\in\{1,2,\ldots,n-1\}\\
n &\mbox{if}& p=q\in\{0,n\}
\end{array}\right.=n \alpha_p \delta_{p,q},
$$
where $$\gamma_{n,k}=\left\{\begin{array}{lcl}
\frac{1}{2} &\mbox{if}& k\in\{0,n\}\\
1 &\mbox{if}& k\in\{1,2,\ldots,n-1\}
\end{array}\right.\ \mbox{and}\ 
\alpha_{n,i}=\left\{\begin{array}{lcl}
\frac{1}{2} &\mbox{if}& i\in\{1,2,\ldots,n-1\}\\
1 &\mbox{if}& i\in\{0,n\}
\end{array}\right.
$$
the Lagrange interpolation polynomial may be written as
$$
L_{n,m}(x,y)=\sum_{i=0}^n\sum_{j=0}^mc_{i,j}T_i(x)T_j(y),
$$
where
$
c_{i,j}=\frac{4}{nm}\gamma_{n,i}\gamma_{m,j}\sum_{k=0}^n\sum_{l=0}^m\gamma_{n,k}\gamma_{m,l}f(x_k,y_l)T_i(x_k)T_j(y_l).
$

\noindent
This polynomial will be called  the Lagrange-Chebyshev interpolating polynomial.
 
As in \cite{1}, the following statements occur:

\begin{teor}\label{cheb2lag11} (Aliasing of Chebyshev polynomials, \cite{1}) For any $n\ge 1$ and $0\le m\le n$ the polynomials
$T_m, T_{n\pm m}, T_{2n\pm m},\ldots$ take the same values on the grid $(\cos{\frac{k\pi}{n}})_{0\le k\le n}.$
\end{teor}

\begin{teor} (Aliasing formula of Chebyshev coefficients)
Let 
$$
f(x,y) = \sum_{i=0}^{\infty}\sum_{j=0}^{\infty}\alpha_{i,j}T_i(x)T_j(y)
$$ 
and let 
$L_{n,m}(x,y)=\sum_{i=0}^n\sum_{j=0}^mc_{i,j}T_i(x)T_j(y)$ be its Lagrange-Chebyshev interpolant.
Then
\begin{equation}\label{cheb2lag10}
c_{i,j}=\sum_{p,q=0}^{\infty}\alpha_{2pn+i,2qm+j}+\sum_{p,q=1}^{\infty}\alpha_{2pn-i,2qm-j}
\end{equation}
\end{teor}

\vspace{0.3cm}\noindent
\textbf{Proof.} Supposing that $(c_{i,j})_{0\le i\le m, 0\le j\le m}$ are given by (\ref{cheb2lag10})
and $\varphi(x,y)= \\\sum_{i=0}^n\sum_{j=0}^mc_{i,j}T_i(x)T_j(y).$
For any $(k,l)\in\{0,1,\ldots,n\}\times\{0,1,\ldots,m\}$
$$
\varphi(x_k,y_l)=\sum_{i=0}^n\sum_{j=0}^mc_{i,j}T_i(x_k)T_j(y_l)=
$$
$$
=\sum_{i=0}^n\sum_{j=0}^m\left(\sum_{p,q=0}^{\infty}\alpha_{2pn+i,2qm+j}+\sum_{p,q=1}^{\infty}\alpha_{2pn-i,2qm-j}\right)T_i(x_k)T_j(y_l).
$$
It is observed that when the indexes $i,j,p,q$ go through their values then $(2pn\pm i,2qm\pm j)$
go through $\mathbb{N}\times\mathbb{N}$ and any two pairs are distinct. Thus, with Theorem \ref{cheb2lag11},
$$
\varphi(x_k,y_l)=\sum_{s=0}^{\infty}\sum_{t=0}^{\infty}\alpha_{s,t}T_s(x_k)T_t(y_l)=f(x_k,y_l).
$$
The unicity of the interpolating polynomial in the set of $(n,m)$ degree polynomials implies $L_{n,m}=\varphi.
\quad \rule{5pt}{5pt}$

A consequence of (\ref{cheb2lag10}) is a relation between $L_{n,m}(x,y)$ and the approximation polynomial
$f_{n,m}(x,y)=\sum_{i=0}^n\sum_{j=0}^m\alpha_{i,j}T_i(x)T_j(y):$
$$
L_{n,m}(x,y)=\sum_{i=0}^n\sum_{j=0}^mc_{i,j}T_i(x)T_j(y)=
$$
$$
=\sum_{i=0}^n\sum_{j=0}^m\left(\alpha_{i,j}+
\sum_{q=1}^{\infty}\alpha_{i,2qm\pm j}+
\sum_{p=1}^{\infty}\alpha_{2pn\pm i,j}+
\sum_{p,q=1}^{\infty}\alpha_{2pn\pm i,2qm\pm j}\right)T_i(x)T_j(y)=
$$
$$
=f_{n,m}(x,y)+\sum_{i=0}^n\sum_{j=m+1}^{\infty}\alpha_{i,j}T_i(x)T_{\mu_j}(y)+
\sum_{i=n+1}^{\infty}\sum_{j=0}^m\alpha_{i,j}T_{\nu_i}(x)T_j(y)+
$$
$$
+\sum_{i=n+1}^{\infty}\sum_{j=m+1}^{\infty}\alpha_{i,j}T_{\nu_i}(x)T_{\mu_j}(y).
$$
Then
\begin{equation}\label{cheb2lag12}
|L_{n,m}-f_{n,m}(x,y)|\le \sum_{i=0}^n\sum_{j=m+1}^{\infty}|a_{i,j}|+
\sum_{i=n+1}^{\infty}\sum_{j=0}^{\infty}|a_{i,j}|.
\end{equation}

Know, it can prove the uniform convergence of the Lagrange-Chebyshev interpolation polynomials:
\begin{teor}
If the function $f$ has second order continuous partial derivatives then $\lim_{n,m\rightarrow\infty}L_{n,m}=f$
uniformly in $[-1,1]^2.$
\end{teor}

\vspace{0.3cm}\noindent
\textbf{Proof.} 
Using (\ref{cheb2lag13}) and (\ref{cheb2lag12}) we obtain 
$$
|f(x,y)-L_{n,m}(x,y)|\le |f(x,y)-f_{n,m}|+|f_{n,m}(x,y)-L_{n,m}(x,y)|\le
$$
$$
\le 2\left(\sum_{i=0}^n\sum_{j=m+1}^{\infty}|a_{i,j}|+
\sum_{i=n+1}^{\infty}\sum_{j=0}^{\infty}|a_{i,j}|\right).
$$
The rest of the proof follows  the proof of Theorem \ref{cheb2lag14}. \rule{5pt}{5pt}

\section{Applications}

1. Integrating $f_{n,m}(x,y)$ on $\Omega$ there is obtained
\begin{equation}\label{che2107}
\intt f(x,y)\mathrm{d}x\mathrm{d}y\approx \intt f_{n,m}(x,y)\mathrm{d}x\mathrm{d}y=4 \sum_{k=0,\  even}^n
\sum_{j=0,\ even}^m \frac{\alpha_{n,m}}{(1-k^2)(1-j^2)}.
\end{equation}

2. Computation of the first order partial derivatives.
Practically, knowing the Chebyshev approximation polynomial
$f_{n,m}(x,y)=\sum_{k=0}^{n}\sum_{j=0}^m\alpha_{k,j}T_k(x)T_j(y),$ and with the
assumption that $\alpha_{k,j}\approx 0$ for $k>n>4,$  and for any $j\in\{0,1,\ldots,m\}$ the first $n$ equations of the system 
(\ref{cheb2111}) will be
$$
\left\{ \begin{array}{lcl}
b_{0,j}-\frac{b_{2,j}}{2} & = & \alpha_{1,j} \\
\frac{1}{2k}(b_{k-1,j}-b_{k+1,j}) & = & \alpha_{k,j},\qquad k\in\{2,3,\ldots,n-2\}\\
\frac{1}{2(n-1)}b_{n-2,j} & = & \alpha_{n-1,j}\\
\frac{1}{2n}b_{n-1,j} & = & \alpha_{n,j}
\end{array}\right.
$$
with the solution
$$
\begin{array}{lcl}
b_{n-1,j} &=& 2n\alpha_{n,j} \\
b_{n-2.j} &=& 2(n-1)\alpha_{n-1,j} \\
b_{k,j} &=& 2(k+1)\alpha_{k+1,j}+b_{k+2,j},\qquad k\in\{n-3,n-4,\ldots,2,1\}\\
b_{0,j} &=& \alpha_{1,j}+\frac{b_{2,j}}{2}
\end{array}.
$$
Then $\frac{\partial f(x,y)}{\partial x}\approx\sum_{k=0}^{n-1}\sum_{j=0}^mb_{k,j}T_k(x)T_j(y).$

The partial derivative $\frac{\partial f(x,y)}{\partial y}$ may be computed similarly.

Due to the truncation of the Chebyshev series the numerical result  is influenced by 
the truncation error as well as by rounding errors. The automatic differentiation \cite{5} is a method
which eliminates the truncation error but it requires a specific computational environment
related to the definition of the \textit{elementary} functions (e.g. \textit{apache commons-math3} v. 3.4).

\section{Examples}

Using a \textit{Scilab} implementation the following results are obtained
\begin{enumerate}
\item
$f(x,y)=\cos{x y}$ \cite{200}, Ch. 11, p. 2.
The matrix of the coefficients is
\begin{verbatim}
0.880725579    0.  - 0.117388011    0.    0.001873213  
     0.        0.    0.             0.    0.           
-0.117388011   0.  - 0.114883808    0.    0.002484444  
     0.        0.    0.             0.    0.           
0.001873213    0.    0.002484444    0.    0.000603385
\end{verbatim}
\normalsize
\noindent
The value of the indicator given by (\ref{che2106}) is $3.97247\cdot 10^{-10}.$  

\noindent
On an equidistant grid of size $50\times 50$ in $[0,1]^2$ the maximum absolute error is 0.000082141.  

\noindent
The integral given by (\ref{che2107}) is 3.784330902, while \textit{Mathematica} gives $4 \mathrm{SinIntegral}[1]\approx 3.78433228147.$
\item
$g(x)=\cos{10 x y^2}+e^{-x^2}$ \cite{200}, Ch. 12, p. 6.
The size of matrix of coefficients is $33\times 43.$

The value of the indicator given by (\ref{che2106}) is 0.  

\noindent
On an equidistant grid of size $50\times 50$ in $[0,1]^2$ the maximum absolute error is  $2.98594\cdot 10^{-13}.$

\noindent
The integral given by (\ref{che2107}) is 4.590369905, which is equal to that given by \textit{Mathematica.}
\end{enumerate}

\section{Conclusions}
There is presented an alternative to the Gaussian elimination method used in \textit{Chebfun} software
in order  to construct an approximation polynomial of a function with two variables.

Because the  discrete Fourier transform is a common tool for the usual mathematical softwares,
this approach has a relative simple implementation, but as a drawback, if the tolerance is the
machine precision then it may require a large amount of memory.

\input cyracc.def
\font\tencyr=wncyr10
\def\cyr{\tencyr\cyracc}

\end{document}